\documentclass{amsart}

\usepackage{graphicx}
\usepackage{amssymb,latexsym}
\usepackage{amsmath,amsfonts,amstext,bbm}

\newtheorem{theorem}{Theorem}
\newtheorem{lemma}{Lemma}

\def\XXint#1#2#3{{\setbox0=\hbox{$#1{#2#3}{\int}$ }
\vcenter{\hbox{$#2#3$ }}\kern-.6\wd0}}

\begin{document}
\title[LEVEL SET ESTIMATES FOR THE DISCRETE FREQUENCY  FUNCTION]
      {Level set estimates for the discrete frequency  function}
      
\author{Faruk Temur}
\address{Department of Mathematics\\
         Izmir Institute of Technology}
\email{faruktemur@illinois.edu}

\date{September 8, 2016}

\maketitle

\begin{abstract}
	We introduce   the  discrete frequency function as a possible new approach to understanding the discrete Hardy-Littlewood maximal function. Considering that the discrete  Hardy-Littlewood maximal function  is given  at each integer by the supremum of averages over intervals of integer length, we define   the discrete frequency function at that  integer as the value at which the supremum is attained.  After verifying that the function is well-defined, we investigate size and smoothness properties of this function.  
\end{abstract}

\section{Introduction}\label{intro}

Let $\mathbb{Z}$ be the set of integers, and let $\mathbb{Z}^+$ denote the set of non-negative integers. Let $f \in l^{1}(\mathbb{Z})$.  For real numbers $a\leq b$,  let $[a,b]$ denote the set of integers $n$ such that $a\leq n\leq b$. We will call $[a,b]$ an interval. We define the average of $f$ over an  interval of radius $r\in \mathbb{Z}^+$ by 
\begin{equation*}
\mathcal{A}_rf(n):= \frac{1}{2r+1}\sum_{k=-r}^{r} |f(n+k)|. 
\end{equation*}
  The discrete Hardy-Littlewood maximal function is defined as   
\begin{equation*}
\mathcal{M}f(n):=\sup_{r\in \mathbb{Z}^+ }\frac{1}{2r+1}\sum_{k=-r}^{r} |f(n+k)| ,
\end{equation*}
thus we have  
\begin{equation*}
\mathcal{M}f(n)=\sup_{r\in \mathbb{Z}^+}\mathcal{A}_rf(n).
\end{equation*}
 Our aim in this work is to study the distribution of the values  $r$ for which $\mathcal{M}f(n)=\mathcal{A}_rf(n)$. More precisely, let 
 \begin{equation*} E_{f,n}:=\{r:\mathcal{M}f(n)=\mathcal{A}_rf(n)\}
 \end{equation*}
  We introduce  the discrete frequency  function as
\begin{equation}\label{f1}
\mathcal{F}f(n):=\inf E_{f,n}.
\end{equation}
 This function is well defined, for the set $E_{f,n}$, which is obviously bounded below, is also non-empty; we will prove this in the next section. Once we have this  we clearly also have $\mathcal{F}f(n)\in E_{f,n}$ since $ E_{f,n}$ is a subset of non-negative integers. We will also prove in the next section that  $E_{f,n}$ is actually finite whenever $f$ is not identically zero.  We will call the  transformation $\mathcal{F}$ the discrete frequency function, for as can be observed  in \cite{ku,t}, the values of $r$ can be used  to decompose the Hardy-Littlewood maximal function in a way that reminds us the decomposition of linear operators  using eigenvalues. And since this decomposition is used to great effect in these works, we find a systematic investigation of this function  very important.  The only investigation of this function that the author could find is \cite{s}, where it is proved that if the frequency function $\mathcal{F}f$ takes only a few values then $f$ must be a sine type function,  although it must be remarked that in that work the functions $\mathcal{M}$ and $\mathcal{F}$ are defined somewhat differently. In this work we will explore aspects of this function quite different from those in   \cite{s}, and we will mainly concentrate on size and smoothness of the frequency function.

We wish to prove level set  estimates for the frequency function, and since we assume $f$ to be  summable, it seems us  most natural to consider level sets obtained by comparing $\mathcal{F}f(n)$ to $|n|$. As being  summable necessitates  decay  at infinity, any such function must have most of its mass on an interval of finite length centered at the origin. This makes the choice of comparison with $|n|$ very natural. But we will also consider  comparisons with various other functions.  We have the following theorems.

\begin{theorem}
Let  $f \in l^1(\mathbb{Z})$. Let $C>1$ be a  real number and let
\begin{equation*}
S_C:=\{n:
\frac{|n|}{2C}\leq\mathcal{F}f(n) \leq \frac{|n|}{C}\}.
\end{equation*}
The set $S_C$ is a finite set.
\end{theorem}

\begin{theorem}
Let  $f \in l^1(\mathbb{Z})$ be a function that is not identically zero. Let $C>1$ be a real number and let
\begin{equation*}
K_{C,N}:=\{n:|n|\leq N, \ 
\mathcal{F}f(n) \leq \frac{|n|}{C}\}.
\end{equation*}
Then 
\begin{equation*}
\lim_{N\rightarrow \infty} \frac{|K_{C,N}|}{N}=0.
\end{equation*}
 We cannot replace $N$ with $ N^{1-\epsilon}$ for any positive $\epsilon$, or even with $ N/ \log^{1+\epsilon}N$.
\end{theorem}

The proof of the second theorem uses a covering lemma that is usually used to prove the classical weak type boundedness result for the maximal function, therefore we suspect that it may be possible to relate this theorem to that result in a relatively short way, although we could not find it.  This would be an important step in understanding both functions. Another important question is having seen that  we cannot replace $N$ in the denominator above with $ N/ \log^{1+\epsilon}N$, whether it is possible to replace it with $ N/ \log N$. This would be another line of inquiry. We point out that Theorem 2 is sharp in another way too. It is not possible to have
\begin{equation*}
\lim_{N\rightarrow \infty} \frac{|K_{C,N,\theta}|\log^{1+\epsilon}N}{N}=0.
\end{equation*}
where $\theta:\mathbb{Z}\rightarrow \mathbb{Z}^+$ is any function satisfying $\theta(n)\leq n/C$, and 
\begin{equation*}
K_{C,N}:=\{n:|n|\leq N, \ 
\mathcal{F}f(n) \leq {\theta(n)}\}.
\end{equation*}
Indeed, we will show this by setting $\theta(n)=0$ for all $n$, which clearly implies all other cases. Therefore relaxing the requirement $\mathcal{F}f(n) \leq |n|/C$ does not give us a better estimate. 

We will also investigate the variational behavior of the frequency function, and show that in fares poorly in this aspect.  We will show that  for any $C>0$ we can find a function $f_C$ such that $\mathcal{F}f_C(1)-\mathcal{F}f_C(0) >C$. By a  more elaborate construction we will also exhibit a function $f$ such that 
\[\sup_{n\in \mathbb{Z}} |\mathcal{F}f(n+1)-\mathcal{F}f(n)|=\infty. \]  

We can define and investigate similar concepts for   the discrete bilinear maximal function as well. 
Let $f,g \in l^1(\mathbb{Z})$. We define for $r\in \mathbb{Z}^+$
\begin{equation*}
\mathcal{B}_r(f,g)(n):= \frac{1}{2r+1}\sum_{k=-r}^{r} |f(n-k)g(n+k)|. 
\end{equation*}
  The bilinear maximal function is defined as   
\begin{equation*}
\mathcal{B}(f,g)(n)=\sup_{r\in \mathbb{Z}^+}\mathcal{B}_r(f,g)(n).
\end{equation*}
  We define the sets
  \begin{equation*}
  E_{f,g,n}:=\{r:\mathcal{B}(f,g)(n)=\mathcal{B}_r(f,g)(n)\}.
  \end{equation*} 
   We introduce  the  function 
\begin{equation}\label{f2}
\mathcal{F}(f,g)(n):=\inf E_{f,g,n}
\end{equation} 
This function is also well defined,  as will be discussed in the next section. It seems reasonable to expect a result analogous to Theorem 2 to hold for this  case as well, but we are not able to prove this. What we are able to show is that there are functions $f,g\in l^1(\mathbb{Z})$ such that for the sets 
\begin{equation*}
K_{C,N}:=\{n:|n|\leq N, \ 
\mathcal{F}(f,g)(n) \leq \frac{|n|}{C}\}.
\end{equation*}
we have 
\begin{equation*}
\lim_{N\rightarrow \infty} \frac{|K_{C,N}|\cdot \log^{1+\epsilon} N}{N}\neq 0.
\end{equation*}

An analogue of the discrete frequency function can be defined for the usual Hardy-Littlewood maximal function that acts on functions on the real line, but since analogues of the sets   $E_{f,n}$ can be empty in that case the definition needs to be more delicate. Furthermore, to prove any kind of level set estimate we need to deal with the issue of Lebesgue measurability. Since these issues make the investigation of that function significantly more complicated, we will carry  that out  in a future  paper.

The rest of the paper proceeds as follows. In the next section we show that both the discrete frequency function and the discrete bilinear frequency function are well defined. In the third section we prove results on the discrete frequency function, and in the fourth we discuss the bilinear discrete frequency function.

\section{Well-Definedness of the Discrete Frequency Functions}
In this section we will show that the discrete frequency functions given by  \eqref{f1} and \eqref{f2} are both  well defined. We start with the function in \eqref{f1}.
As mentioned, this means showing that the set $E_{f,n}$ is non-empty for any summable function $f$ and any integer $n$.
 
 We first note that if the function $f$ is zero everywhere, then the set above obviously is not empty. So we may assume that $f$ is not zero everywhere. In this case for any point $n$ the value $\mathcal{M}f(n)$ is positive. Since we have 
\[ \mathcal{M}f(n)=\sup_{r\in \mathbb{Z}^+}\mathcal{A}_rf(n)\]
 we can find a  non-negative integer $r_1$ such that 
 $\mathcal{M}f(n)-\mathcal{A}_{r_1}f(n)\leq 1 $. Let $d_1$ denote the difference  $\mathcal{M}f(n)-\mathcal{A}_{r_1}f(n)$. Then we can find $r_2\in \mathbb{Z}^+$ such that  $\mathcal{M}f(n)-\mathcal{A}_{r_2}f(n)\leq d_1/2$.
 We thus obtain a sequence $r_1,r_2,r_3,\ldots$ of non-negative integers, and a sequence of  differences $d_1,d_2,d_3,\ldots$ induced by them that satisfy the relation $d_{i+1}\leq d_i/2$.  The set of integers $\{r_i: i\in \mathbb{N}\}$ must be   bounded from above. To see this assume to the contrary that it is not bounded from above. Owing to this assumption we can choose  a subsequence  $r_{i_k}$ of $r_i$  as follows: let $i_1=1$, let $i_2$ be choosen such that $i_2>i_1$, and $r_{i_2}>r_{i_1}$. This is possible for otherwise $r_1$ would be an upper bound for the set $\{r_i: i\in \mathbb{N}\}$. Let $i_3$ be choosen such that $i_3>i_2$, and $r_{i_3}>r_{i_2}$, which is possible for otherwise $\max_{j\leq i_2}r_j$ would be an upper bound for the set $\{r_i: i\in \mathbb{N}\}$.  Thus proceeding  we obtain a subsequence $\{r_{i_k}\}_{k\in \mathbb{N}}$. We clearly have 
 \[\mathcal{M}f(n)= d_{i_k}+ \frac{1}{2r_{i_k}+1}\sum_{j=-r_{i_k}}^{r_{i_k}}|f(n+j)|\leq  d_{i_k}+ \frac{\|f\|_1}{2r_{i_k}+1}.\]
 But  as $k\rightarrow \infty$ the rightmost term converges to zero, while the leftmost term is strictly positive. Therefore  the set  $\{r_{i}:\in \mathbb{N}\}$ must be bounded above. Thus this set actually is finite. Hence for some $r_i$ we must have $\mathcal{M}f(n)=\mathcal{A}_{r_i}f(n)$,
for otherwise $d_i\rightarrow 0$ as $i\rightarrow \infty$ would be impossible.
 
 We now  show that $E_{f,n}$ is finite if $f$ is not identically zero. In this case $\mathcal{M}f(n)$ is strictly positive. If we assume  $E_{f,n}$ to be infinite then we can list  its elements to obtain a sequence $r_1,r_2,r_3,\ldots$ such that $r_1<r_2<r_3<\ldots$. But then 
   \[\mathcal{M}f(n)=  \frac{1}{2r_{i}+1}\sum_{j=-r_{i}}^{r_{i}}|f(n+j)|\leq   \frac{\|f\|_1}{2r_{i}+1}.\]
 Since elements of  $E_{f,n}$ are integers, $r_i\rightarrow \infty$ as $i\rightarrow \infty$. Thus we have a contradiction, and $E_{f,n}$ is finite.
 
 Proof of well-definedness of the bilinear discrete frequency function follows the same  lines. We again wish to prove that the set $E_{f,g,n}$ is not empty. If $\mathcal{B}(f,g)(n)$ is zero then of course $\mathcal{B}_r(f,g)(n)$ is zero for any non-negative $r$, and thus $E_{f,g,n}$ is not empty. So we may assume that $\mathcal{B}(f,g)(n)$ is strictly positive. 
 Since we have 
 \[ \mathcal{B}(f,g)(n)=\sup_{r\in \mathbb{Z}^+}\mathcal{B}_r(f,g)(n)\]
  we can find a  non-negative integer $r_1$ such that 
  $\mathcal{B}(f,g)(n)-\mathcal{B}_{r_1}(f,g)(n)\leq 1 $. Let $d_1$ denote the difference  $\mathcal{B}(f,g)(n)-\mathcal{B}_{r_1}(f,g)(n)$. Then we can find $r_2\in \mathbb{Z}^+$ such that  $\mathcal{B}(f,g)(n)-\mathcal{B}_{r_2}(f,g)(n)\leq d_1/2$.
   we thus obtain a sequence $r_1,r_2,r_3,\ldots$ of non-negative integers, and a sequence of  differences $d_1,d_2,d_3,\ldots$ induced by them that satisfy the relation $d_{i+1}\leq d_i/2$.  The set of integers $\{r_i: i\in \mathbb{N}\}$ must be   bounded from above. To see this assume to the contrary that it is not bounded from above. As before, owing to this assumption we can choose  a subsequence  $r_{i_k}$ of $r_i$  as follows:  $i_1=1$, and for $k\in \mathbb{N}$ we have $i_{k+1}>i_k$ and $r_{i_{k+1}}>r_{i_k}$. Then 
  \[\mathcal{B}(f,g)(n)= d_{i_k}+ \frac{1}{2r_{i_k}+1}\sum_{j=-r_{i_k}}^{r_{i_k}}|f(n-j)g(n+j)|\leq  d_{i_k}+ \frac{\|f\|_1\|g\|_1}{2r_{i_k}+1}.\]
  But  as $k\rightarrow \infty$ the rightmost term converges to zero, while the leftmost term is strictly positive. Therefore we may assume the set  $\{r_i  :i\in \mathbb{N}\}$ to be bounded above. Thus this set actually is finite. Hence  some $r_i$ must be in $E_{f,g,n}$.
 
 We now also prove that if $\mathcal{B}(f,g)(n)$ is not zero then $E_{f,g,n}$ is finite.  If we assume  $E_{f,g,n}$ to be infinite then we can list  its elements to obtain a sequence $r_1,r_2,r_3,\ldots$ such that $r_1<r_2<r_3<\ldots$. Then
     \[\mathcal{B}(f,g)(n)=  \frac{1}{2r_{i}+1}\sum_{j=-r_{i}}^{r_{i}}|f(n-j)g(n+j)|\leq   \frac{\|f\|_1\|g\|_1}{2r_{i}+1}.\]
     Since elements of  $E_{f,g,n}$ are integers, $r_i\rightarrow \infty$ as $i\rightarrow \infty$. Thus we have a contradiction, and $E_{f,g,n}$ is finite.
     
 \section{Proofs of Main Results}
 \subsection{Theorem 1}
 We start with the proof of the first theorem. If  $f$ is identically zero then clearly we have our result. So we will assume $f$ is not identically zero. Assume to the contrary that the set $S_C$ is not finite. Then we have two cases: either positive elements of $S_C$ are infinite, or negative elements of $S_C$ are infinite. We will show the impossibility of the first case, that the second is not possible either can be shown following exactly the same arguments.   Let 
 \[A:=\frac{C+1}{C-1},\ \ \ \ \ B:=\frac{C+1}{C},\ \ \ \ \ D:=\frac{C-1}{C}. \]
 Since $f\in l^1(\mathbb{Z})$ we must have some $m\in \mathbb{N}$ such that 
 \[\sum_{j=-m}^m |f(j)|\geq \frac{\|f\|_1}{2}.\]
  Let $n_1>m$ be a positive element of $S_C$, we can find such an element since we assumed $S_C$ to have infinitely many positive elements.  For the same reason we can find we can find  $n_2\in S_C$ such that $n_2>2An_1$. Proceeding thus we obtain a sequence $\{n_i\}_{i\in \mathbb{N}}$ with $n_{i+1}>2An_i$ for each natural number $i$. Then we observe that
since $n_i\in S_C$ 
\begin{equation*}
\begin{aligned}
\mathcal{M}f(n_i) =\mathcal{A}_{\mathcal{F}f(n_i)}f(n_i) &=\frac{1}{2\mathcal{F}f(n_i)+1}\sum_{j=-\mathcal{F}f(n_i)}^{\mathcal{F}f(n_i)}|f(n_i+j)|\\ &\leq \frac{1}{2\mathcal{F}f(n_i)+1}\sum_{j\in [Dn_i,Bn_i]}|f(j)|.
\end{aligned}
\end{equation*} 
Thus we have 
\[\big(\frac{n_i}{C}+1\big)\mathcal{M}f(n_i) \leq \sum_{j\in [Dn_i,Bn_i]}|f(j)|.\] 
But notice that since $A=B/D$, we have $Dn_{i+1}>2Bn_i$, and therefore the intervals $[Dn_i,Bn_i]$ never intersect. Hence we must have
\begin{equation}\label{t1}
 \sum_{i\in \mathbb{N}} \big(\frac{n_i}{C}+1\big)\mathcal{M}f(n_i) \leq \sum_{i\in \mathbb{N}} \sum_{j\in [Dn_i,Bn_i]}|f(j)|\leq \|f\|_1.
\end{equation}
On the other hand, since 
 $n_i> m$ we must have
\[\mathcal{M}f(n_i) \geq A_{2n_i}f(n_i) = \frac{1}{4n_i+1}\sum_{j=-2n_i}^{2n_i}|f(n_i+j)|=\frac{1}{4n_i+1}\sum_{j=-n_i}^{3n_i}|f(j)|\geq \frac{\|f\|_1}{8n_i+2} .\]
Thus the inequality  \eqref{t1} implies
 \[\frac{\|f\|_1}{C}\sum_{i\in \mathbb{N}} \frac{n_i+C}{8n_i+2}\leq \|f\|_1.\]
Since $f$ is not identically zero, and $C>1$, this implies
\[\frac{1}{C}\sum_{i\in \mathbb{N}} \frac{1}{8}\leq 1,\]
which clearly is not possible.  Thus $S_C$ cannot contain infinitely many  positive elements. 

\subsection{Theorem 2} We now move to the proof of Theorem 2. We will need the following standard covering lemma. By an interval  we mean subsets of integers that contain only consecutive integers, as introduced at the very beginning of this work.

\begin{lemma}
Let $\{B_i\}_{i=1}^m$ be a finite collection of intervals  with  finite length. Let $E$ be a subset of integers   covered by these intervals. Then we can find a disjoint subcollection  $\{B_{i_k}\}_{k=1}^n$ of $\{B_i\}_{i=1}^m$ such that 
\[\sum_{k=1}^n |B_{i_k}|\geq \frac{|E|}{3}.\]
\end{lemma}
This type of lemmas are frequently used to prove boundedness results for maximal functions. For the sake of completeness we will give a proof. Let $B_{i_1}$ be the longest of our intervals. Let $B_{i_2}$ be the longest interval that does not intersect $B_{i_1}$. We choose $B_{i_3}$ to be the longest of the intervals that does not intersect either $B_{i_1}$ or $B_{i_2}$. We proceed thus to obtain a subcollection, which clearly is disjoint. Also observe that any $B_i$ for $1\leq i\leq m$ must intersect  an interval in the subcollection that has at least the same length as itself. For if an interval does not intersect intervals of at least the same length then it must be a member of the collection, which leads to a clear contradiction. Therefore if we consider the collection   $\{3B_{i_k}\}_{k=1}^n$ where $3B_{i_k}$ is the interval obtained by adding a translate of $B_{i_k}$ to its left and  another to its right,  this collection must cover $E$. Therefore 
\[|E|\leq \sum_{k=1}^{n}|3B_{i_k}| = 3 \sum_{k=1}^{n}|B_{i_k}| \]
which clearly implies what we wish.

We can  now start the proof proper. Let $A,B,D$ be defined exactly as in the proof of Theorem 1. Let $K_{C,N}^+$ denote the positive elements of $K_{C,N}$, and $K_{C,N}^-$ denote its negative elements. We will show that 
\[\lim_{N\rightarrow \infty} \frac{|K_{C,N}^+|}{N}= 0,\] 
 and it  will be clear to the reader that the same arguments give this result for $K_{C,N}^-$ as well. Our theorem clearly follows from combining these two results.

 We assume to the contrary that 
\[\lim_{N\rightarrow \infty} \frac{|K_{C,N}^+|}{N}\neq 0.\] 
 This means there exists a small, positive $\epsilon$ such that 
$|K_{C,N_i}^+|/N_i\geq \epsilon$
for a strictly increasing  sequence $\{N_i\}_{i\in\mathbb{N} }$ of natural numbers. So we have 
$|K_{C,N_i}^+|\geq \epsilon N_i$
for such $N_i$. We let $M>10^{10^{10A\epsilon^{-10}}}$  be a natural number such that 
 \[\sum_{j=-M}^M |f(j)|\geq \frac{\|f\|_1}{2}.\]
We now choose a subsequence $\{N_{i_k}\}_{k\in\mathbb{N} }$  of $\{N_i\}_{i\in\mathbb{N} }$ as follows. Let $N_{i_1}$  be such that $N_{i_1}\geq M$, and  let $N_{i_{k+1}}\geq 10A \epsilon^{-1}N_{i_{k}}$ for every $k\geq 1$. Now we fix $k\geq 1$. We have
\[|K_{C,N_{i_{2k}}}^+\setminus K_{C,N_{i_{2k-1}}}^+|\geq \frac{9\epsilon N_{i_{2k}}}{10}\geq 9N_{i_{2k-1}}\]
 Let $n\in K_{C,N_{i_{2k}}}^+\setminus K_{C,N_{i_{2k-1}}}^+ .$  We have 
\[\mathcal{M}f(n)=\frac{1}{2\mathcal{F}f(n)+1}\sum_{j=-\mathcal{F}f(n)}^{\mathcal{F}f(n)}|f(n+j)| \]
but also since no element of the set $K_{C,N_{i_{2k}}}^+\setminus K_{C,N_{i_{2k-1}}}^+$  is in  $[-N_{i_{2k-1}},N_{i_{2k-1}}],$ 
\[\mathcal{M}f(n)\geq \mathcal{A}_{2n}f(n)=\frac{1}{4n+1}\sum_{j=-2n}^{2n}|f(n+j)| =\frac{1}{4n+1}\sum_{j=-n}^{3n}|f(j)|\geq \frac{\|f\|_1 }{8n+2}. \]
Thus combining these two we obtain the fundamental result
\[ \sum_{j=-\mathcal{F}f(n)}^{\mathcal{F}f(n)}|f(n+j)|\geq \frac{2\mathcal{F}f(n)+1 }{8n+2}\|f\|_1. \]
We now consider a covering of $K_{C,N_{i_{2k}}}^+\setminus K_{C,N_{i_{2k-1}}}^+ $ by such $[n-\mathcal{F}f(n),n+\mathcal{F}f(n)]$. By our covering lemma we have a subset $n_1,n_2,\ldots n_{p_k}$ for which the intervals $[n_i-\mathcal{F}f(n_i),n_i+\mathcal{F}f(n_i)], \ 1\leq i\leq p_k$ are disjoint, and
\[\sum_{i=1}^{p_k}2\mathcal{F}f(n_i)+1 \geq\frac{1}{3}|K_{C,N_{i_{2k}}}^+\setminus K_{C,N_{i_{2k-1}}}^+|\geq \frac{9\epsilon N_{i_{2k}}}{30}. \]
We combine this result with the fundamental result above to obtain
\begin{equation*}
\begin{aligned}
\sum_{i=1}^{p_k}\sum_{j=-\mathcal{F}f(n_i)}^{\mathcal{F}f(n_i)}|f(n_i+j)|&\geq \sum_{i=1}^{p_k}\frac{2\mathcal{F}f(n_i)+1 }{8n_i+2}\|f\|_1\\
&\geq \frac{\|f\|_1 }{8N_{i_{2k}}+2}\sum_{i=1}^{p_k}2\mathcal{F}f(n_i)+1\\
&\geq \frac{\|f\|_1 }{8N_{i_{2k}}+2}\frac{9\epsilon N_{i_{2k}}}{30}\\
&\geq \frac{ \epsilon \|f\|_1}{30}
\end{aligned}
\end{equation*}
But since $[n_i-\mathcal{F}f(n_i),n_i+\mathcal{F}f(n_i)]$ are disjoint, we  have
\begin{equation*}
\sum_{j\in[DN_{i_{2k-1}},BN_{i_{2k}}]}|f(j)|\geq \sum_{i=1}^{p_k}\sum_{j=-\mathcal{F}f(n_i)}^{\mathcal{F}f(n_i)}|f(n_i+j)|\geq 
 \frac{ \epsilon \|f\|_1}{30}
\end{equation*}
 Owing to our choice of the subsequence $\{N_{i_k}\}_{k\in \mathbb{N}}$ the intervals $[DN_{i_{2k-1}},BN_{i_{2k}}]$ are disjoint for each natural number $k$, and therefore summing over $k$ we have
\begin{equation*}
\|f\|_1\geq \sum_{k\in \mathbb{N}} \ \ \sum_{j\in[DN_{i_{2k-1}},BN_{i_{2k}}]}|f(j)|\geq 
\sum_{k\in \mathbb{N}} \frac{ \epsilon \|f\|_1}{30}
\end{equation*}
which is a contradiction since $f$ is assumed to be summable and not identically zero.

We now give examples that show the sharpness of the estimate. The following is our  most basic example, and the next two  will improve upon the same ideas. We let for a small, positive $\epsilon$
\[f(n):=\left\{ \begin{split}  &\frac{1}{m^{1+\epsilon}} & \quad &\text{if} \quad n=m^2, \  m\in \mathbb{N},
     \\
      & 0 & \quad &\text{elsewhere} \quad \end{split} \right.\]
Now let $N=M^2$ for  $M>10^{10^{10A\epsilon^{-10}}}$. We have $M^2-(M-1)^2=2M-1$. Let $n$ satisfy $M^2-M^{1-2\epsilon}/4<n<M^2$. We will calculate the maximal function at this point $n$. If we take $r$ to be a natural number satisfying $M^{1-2\epsilon}/2 <r<M^{1-2\epsilon}$, then 
\[\mathcal{A}_rf(n)=\frac{1}{2r+1}\frac{1}{M^{1+\epsilon}}\geq \frac{1}{3M^{1-2\epsilon}M^{1+\epsilon}}=\frac{1}{3M^{2-\epsilon}}\] 
Obviously taking $M^{1-2\epsilon}\leq r<n-(M-1)^2$ cannot give a larger average. We claim that this is not possible for  $r\geq n-(M-1)^2 $ either.   To tackle this case we will use the following observation, which greatly simplifies calculations that otherwise would be very cumbersome. Simply stated our observation is this: as we approach to the origin from the right hand side the function attains nonzero values with increasing frequency, and moreover these nonzero values grow.  In technical terms, we must have average of $f$ over the interval $[(m-2)^2,(m-1)^2-1]$ larger than its average on $[(m-1)^2,m^2-1]$, that is
\begin{equation*}
\begin{aligned}
\frac{1}{2m-3}\sum_{j=(m-2)^2}^{(m-1)^2-1}f(j)&=\frac{1}{(2m-3)(m-2)^{1+\epsilon}} \\ 
 &\geq\frac{1}{(2m-1)(m-1)^{1+\epsilon}}\\&= \frac{1}{2m-1}\sum_{j=(m-1)^2}^{m^2-1}f(j)\\
\end{aligned}
\end{equation*}
Obviously due to this phenomenon  $r$ cannot exceed $n$ too much. Indeed,  a moment's consideration makes it clear that we must have $r<2n$. With such $r$ we must have
\begin{equation*}
\mathcal{A}_rf(n)=\frac{1}{2r+1}\sum_{j=-r}^rf(n+j)\leq \frac{3}{2r+1}\sum_{j=-r}^{-1}f(n+j)\leq \frac{3}{2}\frac{1}{r}\sum_{j=-r}^{-1}f(n+j)
\end{equation*}
Thus at the end we have average over $[n-r,n-1]$ of $f$, and we can write,
\[\frac{1}{r}\sum_{j=-r}^{-1}f(n+j)=\frac{1}{r}\sum_{j=n-r}^{n-1}f(j)\leq2\frac{1}{M^2-n+r}\sum_{j=n-r}^{M^2-1}f(j)\]
 Now we have average over $[n-r,M^2-1]$ of $f$ at the end, and
 we wish to know the greatest value that this average can attain.  Of course if $(m-1)^2<n-r\leq m^2$ for some natural number $m$, taking $r$ so that $n-r=m^2$ makes this  average largest. Then using our observation we conclude that we better take $m=1$. Thus this average is    at most $C_{\epsilon}/M^2$ where $C_{\epsilon}$ is the constant given by
\[C_{\epsilon}=\sum_{m=1}^{\infty}\frac{1}{m^{1+\epsilon}}\leq 2\epsilon^{-1}.\]
 Therefore
$\mathcal{A}_rf(n)\leq3C_{\epsilon}/{M^2} .$
This, given our choice of $M$, is clearly less than $1/3M^{2-\epsilon}$. Thus we must have $\mathcal{F}f(n)\leq |n|/C$. And from amongst $2M-2$ values of $n$ between $(M-1)^2$ and $M^2$, at least $M^{1-2\epsilon}/8$ satisfy this property. If we apply this to each interval $[(M-k-1)^2,(M-k)^2]$  for $k\in [0,M/2]$, we similarly obtain $(M-k)^{1-2\epsilon}/8$ values of $n$ satisfying $\mathcal{F}f(n)\leq |n|/C$. Thus in $[-N,N]$ we have at least 
\[\frac{M}{2}\frac{(M-M/2)^{1-2\epsilon}}{8}\geq  \frac{N^{1-\epsilon}}{50}\]
such elements. Therefore $K_{C,N}$ has at least this cardinality, which makes
\[\lim_{N\rightarrow \infty}\frac{K_{C,N}}{N^{1-\epsilon}}=0  \]
impossible.

We now give our second example. Using exactly the same arguments we can use the function 
 \[f(n):=\left\{ \begin{split}  &\frac{1}{m\log^{1+\epsilon/2} m}& \quad &\text{if} \quad n=m^2, \  m\in \mathbb{N}, \ m\geq 10
      \\
       & 0 & \quad &\text{elsewhere} \quad \end{split} \right.\]
to show that 
\[\lim_{N\rightarrow \infty}\frac{|K_{C,N}|\log^{1+\epsilon}N}{N}=0  \]
is not possible.

Our third example pushes  these ideas to the furthest.  We define for a small positive $\epsilon$
\[f(n):=\left\{ \begin{split}  &\frac{1}{m\log^{1+\epsilon/2} m}& \quad &\text{if}  \quad n=\lceil m\log^{1+\epsilon} m\rceil  , \  m\in \mathbb{N}, \ m\geq 10
      \\
       & 0 & \quad &\text{elsewhere} \quad \end{split} \right. \]
Here  for some real number $x$ the expression $\lceil x \rceil$ denotes the smallest integer that is not less than $x$. We will show  that 
\[\lim_{N\rightarrow \infty}\frac{|K_{C,N,\theta}|\log^{1+\epsilon}N}{N}=0  \]
is not possible for the constant function $\theta(n)=0$.  Let $M>10^{10^{10A\epsilon^{-10}}}$, and let $N=\lceil M\log^{1+\epsilon} M\rceil$. Consider $m\in [M/2,M]$ and values of $n= \lceil m\log^{1+\epsilon} m\rceil$ that correspond to these $m$. For such $m$ we have of course have
\[\mathcal{A}_0f(n)=f(n)=\frac{1}{m\log^{1+\epsilon/2} m}.\]
We will show that $\mathcal{A}_rf(n)$ cannot be larger than this for any $r$. Obviously for $r\gg n$ this is true, indeed a moment's consideration makes it clear that $r< 2n$. For such  $r$ we have
\[\mathcal{A}_rf(n)=\frac{1}{2r+1}\sum_{j=-r}^rf(n+j)\leq \frac{f(n)}{3}+\frac{2}{r}\sum_{j=-r}^{-1}f(n+j).\]
Thus the last term is average over $[n-r,n-1]$, and by the same reasoning as in the first example this average is largest when $r=n-\lceil 10\log^{1+\epsilon} 10\rceil$, for the function attains ever growing nonzero values with ever increasing frequency as we approach to the origin from the right hand side. Therefore
\[\frac{2}{r}\sum_{j=-r}^{1}f(n+j)\leq \frac{4}{n}\sum_{j=10}^{\infty} j\log^{1+\epsilon} j \leq \frac{4C_{\epsilon}}{n}\]
where 
\[\sum_{j=10}^{\infty} j\log^{1+\epsilon/2} j =C_{\epsilon}\leq\frac{2}{\epsilon}.\]
But obviously
\[\frac{4C_{\epsilon}}{n}\leq \frac{16C_{\epsilon}}{M\log^{1+\epsilon} M} \leq \frac{32\epsilon^{-1}}{m\log^{1+\epsilon} m}< \frac{1}{3m\log^{1+\epsilon/2} m}=\frac{f(n)}{3} \]
Therefore $\mathcal{A}_rf(n)< \mathcal{A}_0f(n)$. Thus $K_{C,N,\theta}$ contains at least $M/4$ elements, hence
\[\frac{|K_{C,N,\theta}|\log^{1+\epsilon}N}{N}\geq \frac{M\log^{1+\epsilon}N}{4N}\geq \frac{M\log^{1+\epsilon}M}{8M \log^{1+\epsilon}M}\geq \frac{1}{8}\]
establishing our claim.

\subsection{Variational Results}
For each $C$ positive real number we will show a function $f_C$ such that $\mathcal{F}f_C(1)-\mathcal{F}f_C(0)>C$. Obviously it is enough to find such functions for all $C\in \mathbb{N}, \ C\geq100$. We define for such a $C$
\[f_C(n):=\left\{ \begin{split}  &1& \quad &\text{if}  \quad n=0,
      \\
       & 2C & \quad &\text{if} \quad |n|=3C, 
       \\
        & 0 & \quad &\text{elsewhere} \quad \end{split} \right.\]
Now consider the only reasonable candidates  that may be the value $\mathcal{F}f_C(0)$: the values $0,3C$. We have 
$\mathcal{A}_0f_C(0)=1$ while 
$\mathcal{A}_{3C}f_C(0)=(4C+1)/(6C+1)<1$. Therefore $\mathcal{F}f_C(0)=1$
On the other hand the  only reasonable candidates  that may be the value $\mathcal{F}f_C(1)$ are $1, 3C-1,3C+1$. We have
$\mathcal{A}_1f_C(1)=1/3$ while 
\[\mathcal{A}_{3C-1}f_C(1)=(2C+1)/(6C-1),\ \ \ \mathcal{A}_{3C+1}f_C(1)=(4C+1)/(6C+3).\]
Thus given our large values of $C$ we have $\mathcal{F}f_C(1)=3C+1$ which clearly proves our claim.

We now consider the function 
\[f(n)=\sum_{C=100}^{\infty}2^{-C}f_C(n-4^C)\]
Let $n=4^C$ for some $C \geq 200$. Then obviously only reasonable values for $\mathcal{F}f(n)$ are $0,3C$ or values $r>4^{C-1}$ due to the sparse structure of $f$. We have  again  $\mathcal{A}_0f(n)=2^{-C}$ while 
$\mathcal{A}_{3C}f(n)=2^{-C}(4C+1)/(6C+1)<2^{-C}.$ On the other hand for $r>4^{C-1}$ we have 
\[\mathcal{A}_{r}f(n)=\frac{1}{2r+1}\sum_{j=-r}^rf(n+j)\leq\frac{1}{4^{C-1}}\sum_{j= 100}^{\infty}\frac{4j+1}{2^j}\leq \frac{1}{4^{C-1}}\sum_{j= 1}^{\infty}\frac{1}{(\sqrt{2})^j}\leq \frac{5}{4^{C-1}}\]
which means that $\mathcal{F}f(n)=0$.  Similarly only reasonable values for $\mathcal{F}f(n+1)$ are $0,3C-1,3C+1$ or values $r>4^{C-1}$. Applying exactly the same arguments shows that  $\mathcal{A}_{3C+1}f(n+1)$ is the largest, and therefore $\mathcal{F}f(n+1)=3C+1$. Now since $C$ can be arbitrarily large
\[\sup_{n\in \mathbb{Z}} |\mathcal{F}f(n+1)-\mathcal{F}f(n)|=\infty. \]

\section{The Bilinear Discrete Frequency Function}

In this section we present the  example existence of which we mentioned in the introduction. This example uses the same ideas as in the three examples showing the sharpness of Theorem 2.  We let $f,g$ to be the same function
\[f(n)=g(n):=\left\{ \begin{split}  &\frac{1}{m\log^{1+\epsilon/2} m}& \quad &\text{if}  \quad n=\lceil m\log^{1+\epsilon} m\rceil  , \  m\in \mathbb{N}, \ m\geq 10,
      \\
       & 0 & \quad &\text{elsewhere}. \quad \end{split} \right.\]
 Let $M>10^{10^{10A\epsilon^{-10}}}$, and let $N=\lceil M\log^{1+\epsilon} M\rceil$. Consider $m\in [M/2,M]$ and values of $n= \lceil m\log^{1+\epsilon} m\rceil$ that correspond to these $m$. For such $m$ we have of course have  
\[\mathcal{B}_0(f,g)(n)=f(n)g(n)=\frac{1}{m^2\log^{2+\epsilon} m}\]
We now wish to estimate $\mathcal{B}_r(f,g)(n)$ for  $r$ other than zero. Obviously taking $r>n$ is not reasonable. So  assuming $0<r\leq n$ we have
\[\mathcal{B}_r(f,g)(n)=\frac{1}{2r+1}\sum_{j=-r}^{r}f(n-j)g(n+j)=\frac{1}{2r+1}\sum_{j=-r}^{r}f(n-j)f(n+j).\]
We can write the last sum as
\[\frac{1}{2r+1}\Big[f^2(n)+\sum_{j=-r}^{-1}f(n-j)f(n+j)+\sum_{j=1}^{r}f(n-j)f(n+j)\Big].\]
The last two sums clearly are the same, so we have
\[\frac{1}{2r+1}\Big[f^2(n)+2\sum_{j=1}^{r}f(n-j)f(n+j)\Big].\]
So it is enough to show that 
\[\frac{1}{2r+1}\sum_{j=1}^{r}f(n-j)f(n+j)<\frac{f^2(n)}{3}.\]
We clearly have for $j>1$
\[f(n+j)\leq \frac{1}{m\log^{1+\epsilon/2} m}.\]
Therefore we have
\[\frac{1}{2r+1}\sum_{j=1}^{r}f(n-j)f(n+j) \leq \frac{1}{m\log^{1+\epsilon/2} m}\frac{1}{r}\sum_{j=1}^{r}f(n-j)\]
Thus we again have   the average of $f$ taken over $[n-r,n-1]$ and as explained before this becomes largest when $r=n-\lceil 10\log^{1+\epsilon} 10\rceil$, thus  we have
 \[\frac{1}{m\log^{1+\epsilon/2} m}\frac{1}{r}\sum_{j=1}^{r}f(n-j)\leq \frac{2C_{\epsilon}}{nm\log^{1+\epsilon/2} m}\leq \frac{2C_{\epsilon}}{m^2\log^{2+3\epsilon/2} m} <\frac{f^2(n)}{3} \]
 where 
 \[\sum_{j=10}^{\infty} j\log^{1+\epsilon/2} j =C_{\epsilon}\leq \frac{2}{\epsilon}.\]
Hence we must have at least $M/4$ elements in $K_{C,N}$, and thus
\[\frac{|K_{C,N}|\log^{1+\epsilon}N}{N}\geq \frac{M\log^{1+\epsilon}N}{4N}\geq \frac{M\log^{1+\epsilon}M}{8M \log^{1+\epsilon}M}\geq \frac{1}{8}\]
establishing our claim.

\end{document}